\documentclass[11pt,leqno]{amsart}

\usepackage[utf8x]{inputenc}
\usepackage[T1]{fontenc}
\usepackage[english]{babel}

\usepackage{enumitem}

\usepackage[dvipsnames]{xcolor}
\usepackage{amsthm}
\usepackage[pagebackref, colorlinks=true, linkcolor=MidnightBlue, citecolor=OliveGreen]{hyperref}
\usepackage[capitalise,nameinlink]{cleveref}
\usepackage{amsmath,amssymb}
\usepackage{mathdots}

\usepackage{tikz}

\newtheorem{theorem}{Theorem}[section]
\newtheorem{lemma}[theorem]{Lemma}
\newtheorem{proposition}[theorem]{Proposition}

\newtheorem{question}[theorem]{Question}

\theoremstyle{definition}
\newtheorem{definition}[theorem]{Definition}

\DeclareMathOperator{\Z}{\mathbb{Z}}
\DeclareMathOperator{\N}{\mathbb{N}}
\DeclareMathOperator{\R}{\mathbb{R}}
\DeclareMathOperator{\F}{\mathbb{F}}

\DeclareMathOperator{\Sym}{Sym}

\DeclareMathOperator{\Ball}{B}
\DeclareMathOperator{\depth}{dpt}
\DeclareMathOperator{\diam}{diam}
\DeclareMathOperator{\dist}{d}

\DeclareMathOperator{\id}{id}
\DeclareMathOperator{\Aut}{Aut}
\DeclareMathOperator{\st}{st}
\DeclareMathOperator{\St}{St}
\DeclareMathOperator{\rst}{rist}
\DeclareMathOperator{\Rst}{Rist}
\DeclareMathOperator{\tetr}{tetr}

\DeclareMathOperator{\bin}{bin}

\DeclareMathOperator{\upcom}{U}
\DeclareMathOperator{\leftengel}{LE}

\newcommand{\down}{\hspace{-.3em}\downarrow}

\renewcommand*{\backref}[1]{}
\renewcommand*{\backrefalt}[4]{%
  \ifcase #1 %
    No citations.
  \or
    (cited on page~#2).%
  \else
    (cited on pages~#2).%
  \fi%
}

\makeatletter
\@namedef{subjclassname@2020}{\textup{2020} Mathematics Subject Classification}
\makeatother

\title{On finitely generated Engel branch groups}

\author[J. M. Petschick]{J. Moritz Petschick}
\address{Jan Moritz Petschick: Fakult\"at für Mathematik, Universität Bielefeld, D-33501 Bielefeld, Germany}
\email{jpetschick@math.uni-bielefeld.de}

\date{\today}

\begin{document}
\begin{abstract}
	We construct finitely generated Engel branch groups, answering a question of Fernández-Alcober, Noce and Tracey on the existence of such objects. In particular, the groups constructed are not nilpotent, yielding the second known class of examples of finitely generated non-nilpotent Engel groups following a construction by Golod from 1969. To do so, we exhibit groups acting on rooted trees with growing valency on which word lengths of elements are contracting very quickly under section maps. Our methods apply in principle to a wider class of iterated identities, of which the Engel words are a special case.
\end{abstract}

\maketitle

\section{Introduction} 
\label{sec:introduction}
Given two elements~$g$ and~$h$ of a group~$G$, we set $[g,_0 h] = g$ and furthermore, for~$n \in \N$,
\[
	[g,_{n+1} h] = [[g,_{n} h], h] = [\dots[[g, h], h], \dots h].
\]
The element~$h$ is called \emph{left Engel on~$g$} if there exists some~$n \in \N$ such that $[g,_n h] = \id$. A \emph{left Engel element of~$G$} is an element that is Engel on all~$g \in G$, the set of all such elements is denoted~$\leftengel(G)$, and the group~$G$ is called an \emph{Engel group} if $G = \leftengel(G)$. If, for some $n \in \N$, the identity $[g,_n h]$ is a law in~$G$, one speaks of an $n$-Engel group. The study of Engel elements and Engel groups is a fruitful area of group theory; see~\cite{TT20, Tra11} for some recent surveys. One driving force are the analogues of the family of Burnside problems on periodic groups for the Engel condition. The counterpart of the General Burnside problem may be formulated as the following question, aiming for an analogue of Engel's theorem on finite dimensional Lie algebras:
\begin{enumerate}
	\item Is every finitely generated Engel group nilpotent?
\end{enumerate}
Question~(1) was answered negatively by Golod~\cite{Gol69}. However, his construction remained the only known example of such a group, making it desirable to find new groups that are finitely generated, non-nilpotent, and Engel. In this paper, we construct such groups.

Our examples are branch groups, i.e.\ groups acting on rooted trees whose action resembles the one of the full automorphism group of the tree. The best-known branch groups are the Grigorchuk group and the Gupta--Sidki $p$-groups. One of the remarkable features that these groups share is that they are finitely generated infinite $p$-groups, hence negative answers to the General Burnside Problem. As Golod's construction regarding question~(1) and Golod's counter-example to the General Burnside Problem, see~\cite{Gol65}, use the same methods, it is natural to search for counterexamples regarding question~(1) among branch groups, mirroring the situation for periodic groups. Furthermore, since all branch groups are lawless, every finitely generated Engel branch group is automatically a counter-example to question~(1).

Engel elements in branch groups have been investigated with regards to various aspects, cf.\ \cite{Bar16, FGN19, FNT20, NT18}. Bartholdi~\cite{Bar16} proved that an element of the Grigorchuk group is an Engel element if and only if it satisfies the identity $x^2 = \id$. In particular, since the Grigorchuk group is generated by involutions, the set of left Engel elements is not a subgroup; providing one of the first examples of groups with this behaviour. The first such example, due to Bludov~\cite{Blu05}, was also based on the Grigorchuk group.

In \cite{FGN19}, Fernández-Alcober, Noce and Tracey proved a series of structural results for Engel branch groups; most importantly they established that every such group is a $p$-group for some prime~$p$, strengthening the connexion between question~(1) and the General Burnside Problem.

However, most other results are negative. In fact, aside from the Grigorchuk group and other $2$-groups, no other branch groups with non-trivial Engel elements appear in the literature; furthermore, the fact that $2$-groups such as the Grigorchuk group contain non-trivial Engel elements rests on the fact that for any involution $h$ we find
\(
	[g,_{n+1} h] = [g, h]^{(-2)^{n}}
\)
for sufficiently high $n \in \N$. This method clearly does not extend to other primes. It has been proven that all multi-\textsc{GGS} groups, the Basilica group, the Brunner--Sidki--Vieira group (cf.~\cite{FGN19}), the Hanoi tower group and the group of finitary automorphisms (cf.~\cite{FNT20}) do not contain non-trivial Engel elements. Based on this observation, Fernández-Alcober, Noce and Tracey posed the following questions:
\begin{enumerate}
	\item[(2)] Are there any finitely generated (weakly) branch groups for which the set $\leftengel(G)$ is non-trivial and consists of $p$-elements for an odd prime $p$?
	\item[(3)] Can a finitely generated (weakly) branch group be Engel?
\end{enumerate}
Our examples provide a positive answer to both questions. In fact, we prove the following theorem.

\begin{theorem}\label{thm:Engel groups}
	Let $p$ be a prime. There exists a finitely generated Engel branch $p$-group~$G_p$. In particular, this group is not nilpotent. Furthermore, if the prime~$p$ is odd, the group $G_p$ is contained in a two-generated non-nilpotent Engel group.
\end{theorem}

To achieve this result we define groups where word lengths of elements under sections maps – restrictions to the induced actions on subtrees – are contracting quickly. Contraction properties are used to construct periodic branch groups, cf.\ \cite{Pet23}. The standard definition of a group~$G$ to be contracting is as follows: if there exist $\mu < 1$ and $c \in \R_+$ such that
\[
	\ell(g|_x) \leq \mu\ell(g) + c
\]
holds for all $g \in G$, all first layer section maps $|_x$ and some length function $\ell$ on~$G$. In our case, this definition does not suffice, the problem being that the Engel words $[g,_n h]$ grow very quickly in length for large elements~$h$. We overcome this problem by passing to a more general definition, see \cref{def:contr}, and constructing groups where $\mu$, in a precise sense, is arbitrarily close to~$0$.

This `quick contraction' comes with a price. Our examples act on rooted trees for which the valency of elements is always finite, but unbounded. Since rooted trees with fixed valency -- so called regular rooted trees -- are of particular interest, we prove a weaker statement for groups acting on such trees.

\begin{theorem}\label{thm:Engel p-elements in regular}
	Let~$p$ be a prime. There exists a branch $p$-group acting on a regular rooted tree that is generated by left Engel elements of order~$p$.
\end{theorem}

We suspect that the groups of \cref{thm:Engel p-elements in regular} are not Engel groups. If this is the case, these groups would be the first examples of $p$-groups for an odd prime~$p$ in which the set of left Engel elements is not a subgroup, answering a question posed by Bludov (as cited in~\cite{NT18}) on the existence of such groups.

While the construction of the groups witnessing \cref{thm:Engel groups} and \cref{thm:Engel p-elements in regular} are tailor-made for the Engel property, our general approach allows some generalisation. Using the language of Engel-type iterated identities introduced by Erschler~\cite{Ers15}, we are able to treat being a $p$-group and being an Engel group in a unified way among certain branch groups, cf.\ \cref{thm:strong contraction > engel}, again strengthening the connexion between Engel and periodic groups. Indeed, our methods also apply to a far wider class of iterated identities.

In the hope of stimulating further research, we end this paper with a list of related questions.

\subsection*{Acknowledgements} The author would like to thank Karthika Rajeev for helpful comments on a preliminary version of this text.


\section{Quick contraction and satisfaction of Engel-type identities} 
\label{sec:strong_contraction_and_satisfaction_of_engel_type_identities}

\subsection*{Preliminaries and notation} 
\label{sub:generalities}
Given two integers $n$ and $m \in \Z$ with $m \geq n$, we write $[n, m]$ for the set $\{n, n+1, \dots, m-1, m\}$. To describe certain growth rates, we define the following two integer functions. The \emph{tetration function} $\tetr_n \colon \N \to \N$ with base $n$ is defined by $\tetr_n(0) = 1$ and
\[
	\tetr_n(m+1) = n^{\tetr_n(m)}
\]
for $m \in \N$. Similarly, for two given positive integers~$n > m$, we define the \emph{iterated binomial function} $\bin_{n,m} \colon \N \to \N$ by $\bin_{n,m}(0) = n$ and
\[
	\bin_{n,m}(k) = \binom{\bin_{n,m}(k-1)}{m}
\]
for positive~$k$.

If a group $G$ is generated by a (finite) set $S$, we denote the word length function with respect to $S$ by $\ell_S \colon G \to \N$. We write $\dist_S\colon G \times G \to \N$ for the related distance function, and we write $\Ball_G^S(n)$ for the set of elements of length at most~$n$. The \emph{diameter} of~$G$ with respect to $S$ is defined by $\diam_S(G) = \sup_{g \in G} \ell_S(g)$. Recall that, for two elements~$x$ and~$y$ of a group $G$, we set
\[
	[x,_{n+1} y] = [[x,_{n} y], y]
\]
for $n \in \N$ and $[x_0, y] = x$, where our commutator is conventionally defined by $[x, y] = x^{-1}x^y = x^{-1}y^{-1}xy$.

Let $G$ be a group and let $H$ be a group acting on a set $X$. The wreath product of $G$ and $H$ with respect to the action on $X$ is denoted $G \wr_X H$. Note that in this paper, the group~$H$ (and the set~$X$) will always be finite, whence there is no difference between the restricted and the unrestricted wreath product. The wreath product is, by definition, the semi-direct product $G^X \rtimes H$ with the \emph{top group} $H$ permuting the components of the \emph{base group} $G^X$. For a given $x \in X$, the projection map onto the~$x$\textsuperscript{th}~component is denoted~$\pi_x$. If the set~$X$ is not specified, we refer to the \emph{permutational} wreath product $G \wr H = G^H \rtimes H$ with the top group acting on itself by right multiplication.

Given two integer functions $\mathbf{f}, \mathbf{g}\colon \N \to \R_+$, we define the relation $\prec$ by setting \(\mathbf{f} \prec \mathbf{g}\) if and only if
\[
	\lim_{n \to \infty} \mathbf{f}(n)/\mathbf{g}(n) = 0.
\]
Note that this implies the existence of an integer~$n \in \N$ such that for all $m > n$ we find $\mathbf{f}(m) < \mathbf{g}(m)$.


\subsection*{Automorphisms of rooted trees} 
\label{sub:automorphisms_of_rooted_trees}
Let $\mathcal X = (X_n)_{n \in \N}$ be a sequence of finite sets of cardinality unequal to~$1$. The \emph{rooted (spherically homogeneous) tree} $\mathcal T = \mathcal T_{\mathcal X}$ associated to $\mathcal X$ is the graph with vertex set consisting of all finite strings whose $n$\textsuperscript{th} letter is a member of $X_n$, i.e.\
\[
	\{ x_0x_1 \dots x_n \mid n \in \N, x_i \in X_i \text{ for all } i \in [0, n] \},
\]
including the empty string (or \emph{root}) $\varnothing$ of length $0$, and an edge between strings $x_0 \dots x_n$ and $\tilde x_0 \dots \tilde x_m$ if and only if $|n - m| = 1$ and $x_i = \tilde x_i$ for $i \in [0, \min\{n, m\}]$. The length of a string $u$, denoted $|u|$, is the number of letters in it, and is equal to the length of the unique path in~$\mathcal{T}$ starting at $\varnothing$ and ending at $u$. 
The graph~$\mathcal{T}$ is finite if and only if an empty set occurs in~$\mathcal{X}$.
We will, if not stated otherwise, implicitly assume that $\mathcal T$ is an infinite tree defined by a sequence $\mathcal X = (X_n)_{n \in \N}$ consisting of sets of cardinality at least~$2$.

If the sequence $\mathcal X$ consists of sets of the same cardinality $q \in \N_{>1}$, we call the tree~$\mathcal{T}_{\mathcal X}$ a \emph{$q$-regular rooted tree}.

The group of (graph) automorphisms of a rooted tree~$\mathcal{T}$ is denoted~$\Aut(\mathcal{T})$. For~$g \in \Aut(\mathcal{T})$ and~$u \in \mathcal{T}$, we write~$u.g$ for the image of~$u$ under~$g$.

Given a sequence $\mathcal X$, we denote the $m$-fold shifted sequence $(X_n)_{n \geq m}$ by $\sigma^m \mathcal X$ and the length $m$ prefix-sequence $(X_0, \dots, X_{m-1}, \varnothing, \dots)$ by $\tau^m \mathcal X$.

The set of strings of length $n \in \N$ is called the $n$\textsuperscript{th} layer of $\mathcal T$ and is denoted~$\mathcal L_\mathcal{T}(n)$ or simply $\mathcal{L}(n)$ if it is clear what tree is meant. Let~$G \leq \Aut(\mathcal{T})$ be a group of automorphisms of~$\mathcal{T}$. The stabiliser of a string $u \in \mathcal T$ is denoted $\st_G(u)$. The pointwise stabiliser of~$\mathcal{L}(n)$, the so-called \emph{$n$\textsuperscript{th}~layer~stabiliser}, under the action of a group~$G \leq \Aut(\mathcal T)$ is denoted $\St_G(n)$. The (finite) quotients $G/\St_G(n)$ are called the \emph{congruence quotients of $G$}.

The group $\Aut(\mathcal T)$ has the structure of a wreath product
\[
	\Aut(\mathcal T_\mathcal{X}) = \Aut(\mathcal T_{\sigma \mathcal X}) \wr_{X_0} \Sym(X_0).
\]
The base group of this wreath product is precisely the first layer stabiliser of $\Aut(\mathcal T)$; every component corresponds to automorphisms induced from a subtree isomorphic to $\mathcal{T_{\sigma \mathcal{X}}}$. Iteratively, this yields the identity
\[
	\Aut(\mathcal T) = \Aut(\mathcal T_{\sigma^n \mathcal X}) \wr_{\mathcal L(n)} \Aut(\mathcal T_{\tau^n \mathcal X}).
\]
The tree corresponding to the top group is a finite graph and 
its automorphism group is isomorphic to $\Sym(X_{n-1}) \wr_{X_{n-2}} \dots \wr_{X_0} \Sym(X_0)$. Using the wreath product decompositions of $\Aut(\mathcal{T})$ described above, we define the so-called \emph{section maps}
\[
	|_u \colon \Aut(\mathcal T) \to \Aut(\mathcal T_{\sigma^{|u|}\mathcal X}) \quad\text{ for } u \in \mathcal T
\]
by declaring $g|_u$ to be the projection of~$g$ to the component corresponding to $u \in \mathcal L(n)$. Since~$g$ is not necessarily a member of the base group, this map is in general not a homomorphism. Furthermore, we denote by
\(
	g|^{\mathcal L(n)}
\)
the image of $g$ in the top group $\Aut(\mathcal T_{\tau^n \mathcal X})$, i.e.\ its action on the top part of the tree. Consequently, we may decompose every $g \in \Aut(\mathcal T)$ uniquely as a product
\[
	g = (g|_u)_{u \in \mathcal L(n)} g|^{\mathcal L(n)}.
\]
We record some easily observed identities for section maps. Let $g$ and $h \in \Aut(\mathcal T)$, $u \in \mathcal T$ of length~$n$ and $v \in \mathcal T_{\sigma^n \mathcal X}$, so that the concatenated string $uv$ is naturally an element of $\mathcal T$. Then
\[
	g|_u|_v = g|_{uv} \quad\text{ and }\quad (gh)|_u = g|_u h|_{u.g}.
\]
For a group $G \leq \Aut(\mathcal T)$, the group $\upcom_n(G)$ generated by the sets $G|_u$ for all $u \in \mathcal L(n)$ is called the \emph{$n$\textsuperscript{th} upper companion group of $G$}.

Let $v \in \mathcal L(n)$. The subgroup
\[
	\rst_G(u) = \{ (g|_u)_{u \in \mathcal L(n)} \in \St_G(n) \mid g|_u = \id \text{ for all } v \neq u \}
\]
is called the \emph{rigid stabiliser of $u$}. Evidently all rigid stabilisers of different vertices on the same layer commute. Thus we find
\[
	\prod_{u \in \mathcal L(n)} \rst_G(u) = \langle \rst_G(u) \mid u \in \mathcal L(n) \rangle;
\]
this (normal) subgroup is called the \emph{$n$\textsuperscript{th} rigid layer stabiliser of $G$} and is denoted $\Rst_G(n)$.

A subgroup $G \leq \Aut(\mathcal T)$ is called \emph{spherically transitive} if it acts transitively on every layer $\mathcal L(n)$. In this case all stabilisers $\st_G(u)$ for $u \in \mathcal L(n)$ are conjugate, and the upper companion group $\upcom_n(G)$ is generated by the image of~$G$ under~$|_u$ for any given $u \in \mathcal L(n)$. The group $G$ is called \emph{fractal} if for all $h \in \upcom_n(G)$ and $u \in \mathcal{L}(n)$ there exists $\widehat h \in \St_G(n)$ such that $\widehat h|_u = h$. It is called a \emph{branch group} is it is spherically transitive and every rigid layer stabiliser of $G$ is a finite index subgroup.


\subsection*{Contracting groups} 
\label{sub:_strongly_contracting_groups}
Let $G \leq \Aut(\mathcal T)$ be generated by a set $S \subseteq G$. The \emph{induced generating set for $\upcom_n(G)$} is the finite set of all $n$\textsuperscript{th} level sections of generators $S|_n = \cup_{u \in \mathcal L(n)} S|_u$. To avoid unnecessarily involved notation, for all $m, n \in \N$ we write
\[
	\Ball_G^S(m, n) = \Ball_{\upcom_n(G)}^{S|_n}(m).
\]

\begin{definition}\label{def:contr}
	Let $G \leq \Aut(\mathcal T)$ be a group generated by a set~$S$, let $\mathbf f\colon \N \to \R_{>1}$ be a strictly increasing function and $c \in \R_+$ a positive real number.
	
	The group $G$ is called \emph{contracting with respect to $\mathbf f$ and $c$} if there exists some $n \in \N$ such that for all $m \geq n$, all $g \in G$ and all $u \in \mathcal L(m)$ we have
	\[
		\ell_{S|_m}(g|_u) \leq \frac{\ell_{S|_0}(g)}{\mathbf f(m)} + c.
	\]
\end{definition}

Note that this is a refined version of the usual contraction property; a self-similar group is contracting in the common sense (cf.~\cite{Nek05}) if it is contracting with respect to a function $\mathbf{f}(n) = \mu^{-n}$ for some $\mu \in [0,1)$ and some constant $c$.

It is easy to see that, if given two functions $\mathbf{f}, \mathbf{g}\colon \N \to \N_+$ with $\mathbf{g} \prec \mathbf{f}$ such that $G$ is contracting with respect to~$\mathbf f$ and some constant~$c$, then $G$ is also contracting with respect to $\mathbf g$ and $c$.

Let $G$ be a contracting group (with respect to some function $\mathbf f$) and let $g \in G$. Since~$\mathbf{f}$ is strictly increasing, there exists some $n \in \N$ such that $\mathbf f(m) > \ell_{S}(g)$ for all $m \geq n$, hence we have
\[
	\ell_{S_{n}}(g|_u) \leq \frac{\ell_S(g)}{\mathbf f(n)} + c < 1 + c,
\]
for every $u \in \mathcal L(n)$; thus the length of any section taken at sufficiently long vertices is at most $c$. Following Sidki~\cite{Sid87}, we define the \emph{depth} of an element $g \in G$ as the least integer (or infinity) $\depth(g) \in \N \cup \{\infty\}$ such that $\ell_{S|_{\depth(g)}}(g|_u) \leq c$ for all $u \in \mathcal L(\depth(g))$. The conclusion of the last argument may be rephrased that in a contracting group $G$, the depth $\depth(g) \in \N$ is an integer for all elements $g \in G$.


\subsection*{Satisfaction of Engel-type identities in groups acting on rooted trees} 
\label{sub:satisfaction_of_engel_type_identities}
Let $F_{k+1}$ be the free group on $k + 1$ generators, for some $k \in \N$. We shall denote the generators $x, y_1, \dots, y_k$. Let $w \in F_{k+1}$ an element, i.e.\ a {reduced word in $k + 1$ letters}. Abusing notation, we denote the word map defined by $w$ on a group $G$ by $w \colon G^{k+1} \to G$ by the same symbol. Ultimately, we are interested in Engel and in power words. However, the principal strategy of our proof extends to the more general setting of certain iterated identities.

\begin{definition}(Erschler, \cite[Def~1.1]{Ers15})
	Let $w \in F_{k+1}$ be a word. The sequence $(w_{\circ n})_{n \in \N}$ defined by
	\begin{align*}
		w_{\circ n+1}(x, y_1 \dots, y_k) &= w(w_{\circ n}(x, y_1, \dots, y_k), y_2, \dots, y_k) \quad\text{for } n \in \N \text{ and }\\
		w_{\circ0}(x, y_1, \dots, y_k) &= x.
	\end{align*}
	is called an \emph{Engel-type iterated identity}. A group is called an \emph{Engel-$w$ group} if for all $x, y_1, \dots, y_k \in G$ there exists some $n \in \N$ such that
	\[
		w_{\circ n}(x, y_1, \dots, y_k) = \id.
	\]
	If the $n$\textsuperscript{th} word $w_{\circ n}$ is a law in~$G$, we call~$G$ an \emph{$n$-Engel-$w$ group}.
\end{definition}
In case the commutator word $w = x^{-1}y_1^{-1}xy_1$, we recover the definition of Engel groups, $n$-Engel groups \&c., thus we drop the postfix $w$.

For any word $w \in F_{k+1}$ and for $i \in [1, k]$, we write $a_i(w)$ for the number of times that the generator $y_i$ appears in~$w$, also counting its inverse $y_i^{-1}$, and $a(w)$ for the number of appearances of $x$ or~$x^{-1}$, e.g.\ $a([x, y_1]) = a_1([x, y_1]) = 2$. This allows us to give a (rather crude) upper bound for the length of the images $w_{\circ n}(g, h_1, \dots, h_k)$ under the word map~$w_{\circ n}$ in a group~$G$ generated by some set~$S$, which can be proven easily using induction.

\begin{lemma}\label{lem:upper bound length of iterated word}
	Let $g, h_1, \dots, h_k \in G$ and $n \in \N$. Then
	\[
		\ell_S(w_{\circ n}(g, h_1, \dots, h_k)) \leq a(w)^n \ell_S(g) + \sum_{i = 0}^{n} a(w)^i \sum_{i = 1}^k a_i(w)\ell_S(h_i).
	\]
\end{lemma}


The following lemma will be useful especially when we consider the commutator word, since in this case it allows to restrict to wreath product with cyclic top groups.

\begin{lemma}\label{lem:local checking}
	Let $w \in F_{k+1}$ be a reduced word and let $G \wr_{X} T$ be the wreath product of two groups~$G$ and~$T$ with respect to some action of~$T$ on a set~$X$. Assume that~$T$ is an $n$-Engel-$w$ group. Let
	\[
		g = (g_x)_{x \in X}, h_1 = (h_{1,x})_{x \in X}, \dots, h_k = (h_{k,x})_{x \in X} \in G^X
	\]
	be elements of the base group and let $t_1, \dots, t_k \in T$. For every $x \in X$, denote by $\mathcal O(x)$ the orbit of $x$ under $\langle t_1, \dots, t_k \rangle$.
	Let $m \geq n$. Then the following statements are equivalent:
	\begin{enumerate}
		\item $w_{\circ m}(g, h_{1}t_1, \dots, h_kt_k) = \id$,
		\item $w_{\circ m}((g_x)_{x \in \mathcal O(x')}, (h_{1,x})_{x \in \mathcal O(x')}t_1, \dots, (h_{k,x})_{x \in \mathcal O(x')}t_k) = \id$ for all $x' \in X$.
	\end{enumerate}
	In particular, $G \wr_X T$ is an Engel-$w$ group if and only if $G \wr_{\mathcal{O}(x)} \langle t_1, \dots, t_k \rangle$ is an Engel-$w$ group for all $x \in X$ and all $t_1, \dots, t_k$.
\end{lemma}

\begin{proof}
	Let $R$ be a transversal of the orbits of $X$ under $\langle t_1, \dots, t_k \rangle$. The set $X$ is a disjoint union of the sets $\mathcal O(r)$ for $r \in R$. Order $X$ orbit by orbit and write the tuples within the base group accordingly.
	
	Now collect all top group elements in the expression $w_{\circ m}(g, h_{1}t_1, \dots, h_kt_k)$ to the right. Since $t(g_x)_{x \in X} = (g_{x.t})_{x \in X}$, this does result only in a reordering of the components of the base group elements along some orbit $\mathcal O(r)$. Furthermore, looking at some expression
	\[
		w_{\circ m}((g_x)_{x \in \mathcal O(x')}, (h_{1,x})_{x \in \mathcal O(x')}t_1, \dots, (h_{k,x})_{x \in \mathcal O(x')}t_k) \in G \wr_{\mathcal O(x')}\langle t_1, \dots, t_k \rangle,
	\]
	the same collection process results in the same element of $T$ and in the same base element as in $\mathcal O(x')$-part of the base element above. Thus
	\[
		w_{\circ m}(g, h_{1}t_1, \dots, h_kt_k) = (w_{\circ m}((g_x)_{x \in \mathcal O(r)}, (h_{1,x})_{x \in \mathcal O(r)}t_1, \dots, (h_{k,x})_{x \in \mathcal O(r)}t_k))_{r \in R},
	\]
	where (abusing notation) the tuples on the right hand side are combined into one tuple. This implies the statement.
\end{proof}

Clearly, if any congruence quotient $G/\St_G(n)$ (or any quotient of~$G$ at all) is not an Engel-$w$ group, the group $G$ cannot be an Engel-$w$ group. Consequently, we make the following definition, that may easily be extended to arbitrary groups with a filtration. A group $G \leq \Aut(\mathcal T)$ is called \emph{$\mathbf{f}$-Engel-$w$} for some function $\mathbf f\colon \N \to \N_+$ if $G/\St_G(n)$ is an $\mathbf{f}(n)$-Engel-$w$ group. 

Note that in the profinite topology induced by the congruence filtration, in an $\mathbf{f}$-Engel-$w$ group we have $\lim_{n \to \infty} w_{\circ n}(g, h_1, \dots, h_k) = \id$ for all $g, h_1, \dots, h_k \in G$.

\begin{theorem}\label{thm:strong contraction > engel}
	Let $w \in F_{k+1}$ be a reduced word, let $\mathbf{f}\colon \N \to \N_+$ be some function and let~$G \leq \Aut(\mathcal{T})$ be an $\mathbf{f}$-Engel-$w$ group of automorphisms of some rooted tree~$\mathcal{T}$ that is generated by a finite set~$S$. Assume that~$G$ is contracting with respect to a function~$\mathbf g\colon \N \to \N_+$ such that
	\[
		\max\{ a(w)^{\mathbf{f}}, \mathbf{f}\} \prec \mathbf{g}
	\]
	and a constant~$c \in \R_+$. If there exists some $n \in \N$ such that for all $m \geq n$, every collection of elements $h_1, \dots, h_{k} \in G$ with $\depth(h_j) \leq m$ for all $j \in [1, k]$, every orbit~$\mathcal{O}$ of $\mathcal{L}(m)$ under $T = \langle h_1|^{\mathcal{L}(m)}, \dots, h_{k}|^{\mathcal{L}(m)} \rangle$, and every element $g \in \St_{G}(n)$ such that $\depth(g) \leq m$ the expression
	\[
		w_{\circ i} ((g|_x)_{x \in \mathcal{O}}, (h_1|_x)_{x \in \mathcal{O}} h_1|^{\mathcal{L}(m)}, \dots, (h_k|_x)_{x \in \mathcal{O}} h_k|^{\mathcal{L}(m)}) \in \upcom_m(G) \wr_{\mathcal O} T
	\]
	is trivial for some $i \in \N$, then~$G$ is an Engel-$w$ group.
\end{theorem}

\begin{proof}
	Let $g, h_1, \dots, h_{k} \in G$. We have to show that there exists $n \in \N$ such that $w_{\circ n}(g, h_1, \dots, h_{k}) = \id$. Since $G$ is an $\mathbf{f}$-Engel-$w$ group, we know that
	\[
		w_{\circ \mathbf f(n)}(g, h_1, \dots, h_{k}) \in \St_G(n)
	\]
	for all $n \in \N$. For $u \in \mathcal L(n)$, write $g\down_u$ for $w_{\circ \mathbf f(n)}(g, h_1, \dots, h_{k})|_u$; this is an analogue of the \emph{stabilised sections} defined in \cite{Pet23} for power words $w = x^p$ in order to obtain conditions for periodicity of groups acting on rooted trees. Note that
	\[
		w_{\circ \mathbf f(n)}(g, h_1, \dots, h_{k}) = (g\down_u)_{u \in \mathcal L(n)}.
	\]
	Using \cref{lem:upper bound length of iterated word} and writing $s = \sum_{i=1}^{k} a_i(w) \ell_S(h_i) \in \N$,
	\begin{align*}
		\ell_S(w_{\circ \mathbf f(n)}(g, h_1, \dots, h_{k}))
		&\leq a(w)^{\mathbf f(n)} \ell_S(g) + \sum_{i=0}^{\mathbf{f}(n)} a(w)^i s\\
		&\leq \begin{cases}
			a(w)^{\mathbf{f}(n)}(\ell_S(g) + a(w)s) &\text{ if }a(w) > 1,\\
			\ell_S(g) + (\mathbf{f}(n) + 1) s &\text{ if }a(w) = 1.
		\end{cases}
	\end{align*}
	Since $a(w)^{\mathbf{f}} \prec \mathbf{g}$ (necessary in the first case) and $\mathbf{f} \prec \mathbf{g}$ (for the second case), and since $\ell_S(g), s$ and $a(w)$ are constants, there exists an integer $n_0 \in \N$ such that for all $m > n_0$
	\[
		\ell_S(w_{\circ \mathbf f(m)}(g, h_1, \dots, h_{k}))/\mathbf{g}(m) < 1.
	\]
	Consequently, $\ell_{S|_m}(g\down_u) \leq c$ for all $u \in \mathcal L(m)$. Choose an integer $m$ such that
	\[
		m \geq \max({n_0, \depth h_1, \dots, \depth h_k}),
	\]
	then, by choice, $\ell_{S|_m}(h_i|_u) \leq c$ for all $i \in \{1, \dots, k\}$ and $u \in \mathcal L(m)$. Write $h_i = t_i h_i|^{\mathcal L(m)}$ with $t_i = (t_{1,i}, \dots, t_{|\mathcal L(m)|, i}) \in \Ball(c, m)^{\mathcal L(m)}$.

	For all $j \in \N$ we find
	\[
		w_{\circ m + j}(g, h_1, \dots, h_n) = w_{\circ j}((g\down_u)_{u \in \mathcal L(m)}, t_1 h_1|^{\mathcal L(m)}, \dots, t_k h_k|^{\mathcal L(m)}).
	\]
	In view of~\cref{lem:local checking}, to prove that $w_{\circ i}(g, h_1, \dots, h_k)$ is trivial for some $i \in \N$, we only have to check that for some $i \in \N$
	\[
		w_{\circ i}((g\down_m|_u)_{u \in\mathcal O}, (t_{u,1})_{u \in \mathcal O} h_1|^{\mathcal L(m)}, \dots, (t_{u,k})_{u \in \mathcal O} h_k|^{\mathcal L(m)}) = \id
	\]
	holds for all orbits $\mathcal O$ of $\mathcal L(m)$ under $\langle h_1, \dots, h_k \rangle|^{\mathcal L(m)}$. By our assumption, for sufficiently high~$m$ the last expression is indeed trivial for some $i \in \N$, hence $G$ is an Engel-$w$ group.
\end{proof}



\section{Construction of finitely generated non-nilpotent Engel groups} 
\label{sec:construction_of_finitely_generated_non_nilpotent_engel_groups}

In this section we prove \cref{thm:Engel groups}. Let $p$ be a prime number. We shall construct a tree~$\mathcal T_p$ associated to a sequence of elementary abelian $p$-groups of quickly increasing rank. To start with, we make some basic observations regarding the Cayley graphs of elementary abelian $p$-groups.

If $\mathrm C_p^r$ is an elementary abelian $p$-group of rank $r \in \N$ and $T = \{ e_1, \dots, e_r \}$ a minimal generating set (i.e.\ a basis of $\mathrm C_p^r \cong \F_p^r$), we may define another (except in case $p = 2$ non-minimal) generating set $E = \cup_{i \in [1, r]}\langle e_i \rangle$ consisting of the `coordinate axes' with respect to $T$. It is not difficult to see that the Cayley graphs of $\mathrm C_p^r$ with respect to these generating sets have diameter
\[
	\diam_T(\mathrm C_p^r) = r \lceil (p-1)/2 \rceil
	\quad\text{ and }\quad\diam_E(\mathrm C_p^r) = r,
\]
respectively. We now fix a set of special elements~$F_p(r)$ `far off' the identity element that shall be instrumental later on. For an odd prime $p$, we define $F_p(r)$ to be the set of elements at maximal $T$-distance $\diam_T(\mathrm C_p^r)$ to the identity fulfils
\[
	F_p(r) = \mathrm C_p^r \smallsetminus \Ball_{\mathrm C_p^r}^{T}(\diam_T(C_p^r)-1) = \{ e_1^{\epsilon_1 d} \cdots e_r^{\epsilon_r d} \mid \epsilon_i \in \{-1, 1\} \text{ for } i \in [1, r] \},
\]
where $d = (p-1)/2$ is the $T$-diameter of the rank one group $\mathrm C_p$. Consequently the set~$F_p(r)$ has cardinality $2^r$.

In the case of $p = 2$, the set $F_2(r)$ (if defined as $F_p(r)$ above) has cardinality~$1$, which makes it unfit for our construction. Thus we have to use another definition. Let
\[
	F_2(r) = \Ball_{\mathrm C_2^r}^{T}(r-3) \smallsetminus \Ball_{\mathrm C_2^r}^{T}(r-4),
\]
be the set of elements of `co-length' precisely~$3$. For $r > 3$, this set has cardinality $\binom{r}{3}$ and consists of the elements of the form $e_1^{\epsilon_1} \dots e_r^{\epsilon_r}$ with $\epsilon_{i_1} = \epsilon_{i_2} = \epsilon_{i_3} = 0$ for precisely three indices, and $\epsilon_i = 1$ otherwise.

We now define the trees on which our groups act. Let $X_{0,p}$ be the group of order~$p$, generated by the element~$e_{\varnothing,p}$, for odd~$p$ and $X_{0,2} = \langle e_{1,2}, \dots, e_{5,2} \rangle \cong \mathrm C_2^5$ an elementary abelian $2$-group of rank~$5$. Assume that the set~$X_{n,p}$ is already defined for some~$n \in \N$. Then $X_{n+1,p}$ is the elementary abelian $p$-group (minimally) generated by the symbols $e_{f,p}$ for $f \in F_p(n)$, i.e.\ is isomorphic to $\mathrm C_p^{2^{\operatorname{rk}(X_{n,p})}}$ or, for $p = 2$, the elementary abelian $2$-group of rank~$\binom{\operatorname{rk}(X_{n,2})}{3}$. Consequently
\[
	\operatorname{rk}(X_{n,p}) = \begin{cases}
		\tetr_2(n) &\text{ if } p \neq 2,\\
		\bin_{5,3}(n) &\text{ if } p = 2.
	\end{cases}
\]
Set $\mathcal X_p = (X_{n,p})_{n \in \N}$ and $\mathcal T_p = \mathcal T_{\mathcal X_p}$. An element $x \in X_{n,p}$, for some $n \in \N$, acts by right-multiplication on $X_{n,p}$. This action induces an action on the tree $\mathcal T_{\sigma^n \mathcal X_p}$ by setting $(x_n \dots x_{n+t}).x = (x_n \cdot x)x_{n+1}\dots x_{n+t}$ for all $x_n \dots x_{n+t} \in \mathcal T_{\sigma^n \mathcal X_p}$. Consequently we identify every vertex of layer $n$ with an automorphism of the corresponding tree~$\mathcal T_{\sigma^n \mathcal X_p}$.

We now define a sequence $(b_{n,p})_{n \in \N}$ of elements of $\St_{\Aut(\mathcal T_{\sigma^n \mathcal X_p})}(1)$ by
\[
	b_{n,p}|_x = \begin{cases}
		b_{n+1,p} &\text{ if } x = \id,\\
		e_x &\text{ if } x \in F_p(n),\\
		\id &\text{ otherwise.}
	\end{cases}
\]
It is not difficult to see that $b_{n,p}$ is of order~$p$: its $p$\textsuperscript{th} power fixes the first layer, obviously acts trivial on all subtrees except the subtree below the vertex $\id$, where it acts by the $p$\textsuperscript{th} power of $b_{n+1,p}$. Repeating the argument, we see that~$b_{n,p}^p$ fixes every vertex.

For an odd prime~$p$, let~$\mathcal{G}_{p}$ be the group generated by the set $E_p = \langle b_{0,p} \rangle \cup \langle e_{\varnothing,p} \rangle$. This group is minimally generated by $\{b_{0,p}, e_{\varnothing,p}\}$). Analogously, let $\mathcal G_2$ be the group generated by $E_2 = \langle b_{0,2} \rangle \cup \langle e_{1,2} \rangle \cup \dots \cup \langle e_{5,2} \rangle$.

At this point one may notice that the diverging definitions for $p = 2$ are necessary, since the construction used for odd primes yields a group generated by two elements of order~$p$; thus for $p=2$ we would just obtain a dihedral group.

From here on, we drop the subscript~$p$ from all objects defined above to declutter the notation; when useful or necessary for the clarity of the exposition, we will reintroduce the subscript from time to time.

As a shorthand, denote the upper companion groups $\upcom_n(\mathcal{G})$ by $\mathcal{G}_{n}$. We shall prove that $\mathcal{G}$ is a non-nilpotent Engel group, and that $\mathcal{G}_{n}$ is an Engel branch group for sufficiently high values of~$n$.

\begin{lemma}\label{lem:basic}
	The group $\mathcal G_{n}$ is spherically transitive and fractal for all non-negative integers~$n \in \N$.
\end{lemma}

\begin{proof}
	Since the upper companion groups of spherically transitive fractal groups are again spherically transitive and fractal, we may restrict our considerations to the full group~$\mathcal{G}$.
	
	The upper companion groups of $\mathcal G$ are generated by the induced generating sets $E|_n = \langle b_n \rangle \cup \bigcup_{f \in F(n-1)} \langle e_f \rangle$. To prove that $\mathcal G$ is fractal, it is enough to show that $E|_n \subseteq \st_{\mathcal G}(u)|_u$ for every $u \in \mathcal L(n)$. Now $E|_n$ consists of elements of order at most~$p$, so if some non-trivial power of a generator is contained in $\st_{\mathcal G}(u)|_u$, the generator itself is in the stabiliser as well. Thus it is enough to prove that $b_n$ and $e_f$ for all $f \in F(n-1)$ are contained in $\st_{\mathcal G}(u)|_u$.
	
	Let $u \in \mathcal L(1) = X_0$. Then, viewing $u \in X_0$ as an element of $\mathcal G$ and, recalling that $b_0 \in \St_{\mathcal G}(1) \leq \st_{\mathcal G}(u)$, we see that
	\[
		(b_0)^{u}|_u = b_0|_{uu^{-1}} = b_1 \quad\text{and}\quad (b_0)^{f^{-1}u}|_u = b_0|_{uu^{-1}f} = e_f
	\]
	for all $f \in F(0)$. Given $u = v x_{n} \in \mathcal L(n)$ with $x_n \in X_{n-1}$ and $v \in \mathcal L(n-1)$ we may suppose inductively that there are $\widehat b$ and $\widehat x \in \st_{\mathcal G}(v)$ such that $\widehat{b}|_v = b_{n-1}$ and $\widehat{x}|_v = x$ for all $x \in X_{n-1}$. Clearly $b_{n-1}$ stabilises $u$, and, as above, we find
	\[
		\widehat{b}^{\widehat{x_n}}|_u = (\widehat{b}|_v)^{\widehat{x_n}|_v}|_{x_n} = (b_{n-1})^{x_n}|_{x_n} = b_n
	\] and $\widehat{b}^{\widehat{f^{-1}x_n}}|_u = e_f$ for $f \in F(n-1)$. Thus $\mathcal G$ is fractal.
	
	Now $X_n$ acts transitively on itself and $X_n \leq \st_{\mathcal G}(u)|_u$ for all $u \in \mathcal L(n-1)$, from which we deduce that $\mathcal G$ is spherically transitive.
\end{proof}

Besides the generating set~$E$ and the induced generating sets $E|_n$ we shall use another sequence of generating sets for $\mathcal G$ and its upper companion groups. For all $n \in \N$, set $S|_n = \langle b_n \rangle^{X_n} \cup X_n$; it is easy to observe that the sets $S|_n$ are indeed the induced generating sets for $S = S|_0$.

\begin{lemma}\label{lem:S to E}
	Let $n \in \N$ be an integer, let $g \in \mathcal G_{n}$ an element, and let $x \in X_{n}$ be a vertex. Then $\ell_{E|_{n+1}}(g|_x) \leq \ell_{S|_n}(g)$.
\end{lemma}

\begin{proof}
	Every reduced word in $S|_n$ evaluating to $g$ may easily be transformed into a word of the form $(b_n^{k_1})^{x_1} \dots (b_n^{k_{t-1}})^{x_{t-1}} x_t$, with $x_i \in X_n$ (possibly being trivial) and $k_i \in [1, p-1]$ for $i \in [1, t]$, and the same $S|_n$-length. Since $b_n \in \St_{\mathcal G_{n}}(1)$, we have
	\[
		g|_x = (b_n^{k_1})^{x_1}|_x \dots (b_n^{k_{t-1}})^{x_{t-1}}|_x x_t|_x = (b_n^{k_1})|_{xx_1^{-1}} \dots (b_n^{k_{t-1}})|_{xx_{t-1}^{-1}} x_t|_x.
	\]
	First layer sections of powers of $b_n$ are at most of $E_{n+1}$-length $1$, while all proper sections of~$x_t$ are all trivial. Thus $\ell_{E_{n+1}}(g|_x) \leq t-1 \leq \ell_{S|_n}(g)$.
\end{proof}

To describe the speed of contraction of~$\mathcal G$, we define a function $\mathbf{d}\colon \N \to \N_+$ by
\[
	\mathbf{d}(n) = \begin{cases}
		\tetr_2(n) &\text{ if }p \neq 2,\\
		\bin_{5,3}(n)-3 &\text{ if } p = 2.
	\end{cases}
\]
Note that $\mathbf{d}(n) = \min_{ f \in F(n) } \ell_{{E|_n}}(f) = \ell_{{E|_n}}(f)$ for any $f \in F(n)$. Building on this function $\mathbf{d}$, we furthermore define $\mathbf{g}: \N \to \N$ by $\mathbf{g}(0) = \mathbf{g}(1) = 1$ and, for all $n > 1$, by
\[
	\mathbf{g}(n) = 
		\prod_{{[0, n-1] \ni i \equiv_2 n}} \mathbf{d}(i).
\]
Evidently both functions implicitly depend on $p$. For odd $p$, the functions $\mathbf{d}$ and $\mathbf{g}$ are evidently asymptotically exceeding any exponential function, and for $p = 2$ we find
\[
	(n \mapsto 2^{(2^{n})}) \prec \mathbf{d} \prec \mathbf{g}.
\]

\begin{proposition}\label{prop:strong contr}
	The group~$\mathcal{G}$ is contracting with respect to~$\mathbf{g}$ and~$1$.
\end{proposition}

\begin{proof}
	Let $n \in \N$. We claim that $\ell_{E|_{n+2}}(g|_u) \leq \left\lceil\frac {\ell_{E|_n}(g)} {\mathbf{d}(n)} \right\rceil$ for all $g \in \mathcal G_{n}$ and $u \in X_{n}X_{n+1}$. Since $\ell_{E|_1}(g|_u) \leq \ell_{S|_0(g)} \leq \ell_{E}(g)$ for all $u \in \mathcal L(1)$ by \cref{lem:S to E} and the fact that $E|_n \subset S|_n$, this shows that
	\[
		\Ball_{\mathcal G}^{E}(\mathbf{g}(n))|_u \subseteq \Ball_{\mathcal G}^{E}(1, n)
	\]
	for all $u \in \mathcal L(n)$ by iterating the claimed inequality. Factorising any element $g \in \mathcal G$ into elements of length at most $\mathbf{g}(n)$ we deduce
	\[
		\ell_{E|_{n}}(g|_u) \leq \frac {\ell_{E}(g)}{\mathbf{g}(n)} + 1,
	\]
	which proves the proposition.
		
	To prove the claim it is enough to show that every $g \in \Ball^{E}_{\mathcal G}(\mathbf{d}(n), n)$ satisfies
	\[
		\ell_{E|_{n+2}}(g|_u) \leq 1
	\]
	for all $u \in X_{n}X_{n+1}$, employing the same argument as above. Thus let $g \in \Ball^{E}_{\mathcal G}(\mathbf{d}(n), n)$. There is nothing to prove for $n = 0$ and for $n = 1$ since, again, the statement follows from \cref{lem:S to E}. Assume $n > 1$. Starting with any $S|_n$-word representing $g$, we collect all generators $x \in X_n$ to the right, and obtain a $S|_{n}$-word representing $g$ of at most the same length,
	\begin{equation*}\label{eq:sec prod}
		g = (b_n^{k_1})^{x_1} (b_n^{k_2})^{x_2} \dots (b_n^{x_{t-1}})^{x_{t-1}} x_{t}.\tag{$\ast$}
	\end{equation*}
	Again we have $\ell_{E|_n}(g) \geq \ell_{S|_n}(g) \in \{t-1, t\}$ by $E|_n \subseteq S|_n$. Since for all $x \in X_{n+1}$
	\begin{equation*}
		g|_x = (b_n^{k_1})^{x_1}|_x \dots (b_n^{k_{t-1}})^{x_{t-1}}|_x,
	\end{equation*}
	we have to consider the sections of the elements in $\langle b_n \rangle^{X_n} \subset S|_n$.
	For all $k \in \{1, \dots, p-1\}$ and $y \in X_n$
	\begin{equation*}\label{eq:syl sects}
		(b_n^{k})^{y}|_x = b_n^{k}|_{xy^{-1}} = \begin{cases}
			b_{n+1}^k &\text{ if }x = y\\
			e_f^k &\text{ if }x = fy \text{ for } f \in F(n)\\
			\id &\text{ otherwise}.
		\end{cases}\tag{$\dagger$}
	\end{equation*}
	Thus all factors in the product \eqref{eq:sec prod} that are not of the form $(b_n^k)^{y}$ or $(b_n^k)^{fy}$ for some $f \in F(n)$ do not contribute to $g|_x$. Since we want to maximise the length of $g|_x$, we may assume that all factors are of the kind $(b_n^k)^y$ or $(b_n^k)^{fy}$.
	
	If all factors evaluate to elements in $\langle b_{n+1} \rangle$ or to elements in $X_{n+1}$, respectively, the $S_{n+1}$-length of $g|_x$ is at most $1$. Thus assume that at least one factor of both types appears. Without loss of generality, some $b_{n+1}^i$ appears first and some $e_f^j \in X_{n+1}$ appears second, for $i, j \in [1, p-1]$ and $f \in F(n)$, hence
	\[
		g = (b_n^{i})^{x} (b_n^{j})^{f^{-1}x} \tilde g = x^{-1} b_n^i x f b_n^j \tilde g
	\]
	for some $\tilde g \in \mathcal G_{n}$ such that the second expression is ${E|_{n}}$-reduced. Thus, keeping in mind $\ell_{E_n}(f) \leq \ell_{E_n}(xf) + \ell_{E_n}(x)$, we have
	\begin{align*}
		\mathbf{d}(n) &= \ell_{E|_n}(g) = \ell_{E|_n}(x) + \ell_{E|_n}(b_n^i) + \ell_{E|_n}(xf) + \ell_{E|_n}(b_n^j) + \ell_{E|_n}(\tilde g)\\
		&\geq \ell_{E|_n}(f) + 2 + \ell_{E|_n}(\tilde g)
		\geq \ell_{E|_n}(f) + 2 = \begin{cases}
			\tetr_2(n) + 2 &\text{ if } p \neq 2,\\
			\bin_{5,3}(n) - 1  &\text{ if } p = 2
		\end{cases}\\
		&\geq \mathbf{d}(n) + 2,
	\end{align*}
	a contradiction. Hence $\ell_{S|_{n+1}}(g|_x) \leq 1$ and, by \cref{lem:S to E}, $\ell_{E|_{n+2}}(g|_u) \leq 1$ for all $u \in X_{n}X_{n+1}$.
\end{proof}

\begin{lemma}\label{lem:max orbit}
	The maximal size of an orbit of a vertex $u \in \mathcal L(n)$ under the action of an element $g \in \mathcal G$ is $p^{n}$.
\end{lemma}

\begin{proof}
	This is clearly true for $n = 0$. Assume the statement is true for some $n \in \N$. Hence $g^{p^n} \in \St_{\mathcal G}(n)$. Thus $g^{p^n} = (g^{p^n}|_u)_{u \in \mathcal L(n)}$. But $(g^{p^n}|_u)^p \in \St_{\mathcal G_{n}}(1)$, since $\mathcal G_{n}$ acts on $X_n$ as an elementary abelian $p$-group, hence $g^{p^{n+1}} \in \St_{\mathcal G}(n+1)$.
\end{proof}

\begin{proposition}\label{prop:p group}
	The group $\mathcal G$ is a $p$-group.
\end{proposition}

Note that \cref{prop:p group} follows directly from \cref{prop:engel} and \cite[Theorem~B]{FNT20}. Alternatively, it may be deduced with some care from the proof of \cite[Theorem~A]{Pet23}. However, a proof using \cref{thm:strong contraction > engel} is strait-forward.

\begin{proof}
	In the terminology of \cref{thm:strong contraction > engel}, we set $w = x^p$. Thus $k = 0$ and $a(w) = p$. \cref{lem:max orbit} states that $\mathcal G$ is a $\mathbf{f}$-Engel-$w$ group for $\mathbf{f}(n) = p^n$. Clearly $\mathbf{f} \prec p^{\mathbf{f}} \prec \mathbf{g}$. Since the orbits under $k$-generated -- i.e.\ trivial -- subgroups are of size~$1$, it remains to show that for some $n \in \N$ and all $m \geq n$ the set $\Ball_{\mathcal G}^{E}(1, n)$ consists of $p$-elements. But this is clearly true for all $n$.
\end{proof}

The following lemma represents the crux of our argument. It says, heuristically, that for an element $g$ of sufficiently small word length, sections that are `close' to each other commute. Since we may, using \cref{thm:strong contraction > engel}, restrict our attention to words of bounded length, this will allow us to prove that the group $\mathcal G$ is an Engel group.

\begin{lemma}[Separation]\label{lem:separation}
	Let $n \in \N$ be an integer, let $t \leq \mathbf{d}(n)/2$, and let $x \in X_n$. Then
	\[
		\Ball_{\mathcal G}^{E}(t, n)|_x \subseteq
		\begin{cases}
			\langle b_{n+1} \rangle &\text{ if }\ell_{E|_{n}}(x) \leq t,\\
			X_{n+1} &\text{ if }\ell_{E|_{n}}(x) \geq \mathbf{d}(n) - t,\\
			\{ \id \} &\text{ otherwise.}
		\end{cases}
	\]
	In particular, if $g, h$ are in $\Ball_{\mathcal{G}}^{E}(t, n)$ and $y,z \in X_n$ such that $[g|_y, h|_z] \neq \id$, then $\dist_{{E|_n} \cup X_n}(y,z) \geq \mathbf{d}(n) - 2t$.
\end{lemma}

\begin{proof}
	The idea used here is similar to the proof of \cref{prop:strong contr}. Let $g \in \Ball_{\mathcal G}^{E}(t, n)$ be represented by a reduced $S|_n$-word $(b_n^{k_1})^{x_1} \dots (b_n^{k_{s-1}})^{x_{s-1}} x_s$ with $k_i \in [1, p-1]$ and $x_i \in X_{n}$ for $i \in [1, s]$. Since $\ell_{E|_n}(g) \leq t$, we have $\ell_{E|_n}(x_i) \leq t$ for all $i \in [1,s]$. Looking at \eqref{eq:syl sects} in the proof of \cref{prop:strong contr}, we find that the first case is only possible if $x \in \{x_i \mid i \in [1, s-1]\}$, hence $\ell_{E|_n}(x) \leq t$, and the second case is only possible if $x = fx_i$ for some $i \in [1, s-1]$ and $f \in F(n)$, hence if $\ell_{E|_n}(x) \geq \mathbf{d}(n) - t$.	
\end{proof}

\begin{lemma}\label{lem:vanishing commutators}
	Let $g \in \Ball_{\mathcal G}^{E}(t, n)$ for $t \leq \mathbf{d}(n)/4 - 1$. Then
	\[
		[b_n, g, b_n] = \id.
	\]
	In particular, this is true for all $n \in \N$ and $g \in \Ball_{\mathcal G}^{E}(1, n+3)$.
\end{lemma}

\begin{proof}
	Note that $[b_n, g]$ and $b_n$ are contained in $\St_{\mathcal G_{n}}(1)$. Thus, we may compute their commutator section-wise, and obtain
	\[
		[b_n, g, b_n] = ([[b_n, g]|_x, b_n|_x])_{x \in X_n}.
	\]
	Since $\ell_{E|_n}([b_n, g]) \leq 2(\ell_{E|_n}(g) + 1) \leq \mathbf{d}(n)/2$, the first statement follows easily from \cref{lem:separation}. The additional claim is a consequence of the readily observed fact that $\mathbf{d}(3) \geq 16$.
\end{proof}

\begin{proposition}\label{prop:branch}
	The group~$\mathcal G_{n}$ is a branch group for~$n \geq 2$.
\end{proposition}

\begin{proof}
	Since the upper companion groups of branch groups are branch groups themselves, it is enough to prove the claim for~$\mathcal G_{2}$. Note that by \cref{lem:basic}, the group~$\mathcal G_{2}$ is spherically transitive. Furthermore, its upper companion groups are the same as ones of~$\mathcal{G}$, up to a shift of index, i.e.\ $\upcom_n(\mathcal{G}_2) = \upcom_{n+2}(\mathcal{G}) = \mathcal{G}_{n+2}$.
	
	Let $n \in \N$ and $u \in \mathcal L_{\mathcal T_{\sigma^2 \mathcal{X}}}(n)$. We prove that $\rst_{\mathcal G_{2}}(u) \geq \gamma_3(\mathcal G_{n+2})$. For $n = 0$ this is clear. We proceed inductively, assuming the statement for some $n \in \N$. Consider thus the group~$\mathcal{G}_{n+3}$ is generated by the set $\{ b_{n+3} \} \cup \{e_f \mid f \in F(n+2)\}$. The derived subgroup $\mathcal{G}_{n+3}'$ is normally generated (within $\mathcal{G}_{n+3}$) by $\{ [b_{n+3}, e_f] \mid f \in F(n+2) \}$, and the next term of the lower central series $\gamma_3(\mathcal G_{n+3})$ is normally generated by
	\[
		\{ [b_{n+3}, e_f, b_{n+3}] \mid f \in F(n+2) \} \cup \{ [b_{n+3}, e_f, e_{f'}] \mid f, f' \in F(n+2) \}.
	\]
	By \cref{lem:vanishing commutators}, the first set may be dropped. This is the only point where we use the fact that we are dealing with $\mathcal{G}_{2}$ rather than $\mathcal{G}$. By induction, we may assume that the elements $[b_{n+2}, b_{n+2}^{f^{-1}}, b_{n+2}^{f'^{-1}}] \in \gamma_3(\mathcal G_{n+3})$ are contained in $\rst_{\mathcal G_{2}}(v)$, where $vx = u$ for a vertex~$v$ of length~$n$ and a letter $x \in X_{n+3}$. Since $\mathcal{G}_{2}$ is spherically transitive, we may assume $x = \id_{X_{n+2}}$. We shall evaluate these elements below $v$ and find all first layer sections aside from the one at $\id_{X_{n+3}}$ trivial. Note first that all sections of an element of the form $b_{n+2}$ or $b_{n+2}^f$ -- except the unique section equal to $b_{n+3}$ -- commute. Thus all sections of $[b_{n+2}, b_{n+2}^{f^{-1}}, b_{n+2}^{f'^{-1}}]$ that do not contain $b_{n+3}$ are trivial. Hence for all $x \in X_{n+3}$
	\[
		[b_{n+2}, b_{n+2}^{f^{-1}}, b_{n+2}^{f'^{-1}}]|_x = 
		\begin{cases}
			[b_{n+3}, e_f, e_{f'}] &\text{ if }x = \id_{X_{n+3}},\\
			[e_f, b_{n+3}, b_{n+2}|_{f^{-1}f'}] &\text{ if }x = f^{-1},\\
			[e_{f'}, b_{n+2}|_{f'^{-1}f}, b_{n+3}] &\text{ if }x = f'^{-1},\\
			\id &\text{ otherwise}.
		\end{cases}
	\]
	If $f \neq f'$, the elements $b_{n+2}|_{f'^{-1}f}$ and $b_{n+2}|_{f^{-1}f'}$ are trivial. If $f = f'$, the section $[b_{n+2}, b_{n+2}^{f^{-1}}, b_{n+2}^{f'^{-1}}]|_{f^{-1}} = [e_f, b_{n+3}, b_{n+3}] = \id$. Now let $g \in \mathcal G_{n+3}$. By \cref{lem:basic}, there exists $\widehat g \in \St_{\mathcal G_{2}}(u)$ such that $\widehat g|_u = g$. Let $h \in \gamma_3(\mathcal G_{n+3})$. There exist $k \in \N$, $g_1, \dots, g_k \in \mathcal G_{n+3}$ and $f_1, f'_1, \dots, f_k, f'_k \in F(n)$ such that
	\begin{align*}
		h &= [b_{n+3}, e_{f_1}, e_{f'_1}]^{g_1} \dots [b_{n+3}, e_{f_k}, e_{f'_k}]^{g_k}\\
		&= ([b_{n+2}, b_{n+2}^{f_k^{-1}}, b_{n+2}^{{f'}_k^{-1}}]^{\widehat g_1})|_{\id}  \dots ([b_{n+2}, b_{n+2}^{f_k^{-1}}, b_{n+2}^{{f'}_k^{-1}}]^{\widehat g_k})|_{\id}.
	\end{align*}
	Clearly all other sections of the factors in the second product are trivial. Thus $h \in \rst_{\mathcal G_{2}}(u)$.
	
	It remains to show that $\Rst_{\mathcal G_{2}}(n)$ is of finite index in $\mathcal{G}_2$. We have
	\[
		\gamma_3(\mathcal G_{n+2})^{\mathcal L(n)} \leq \Rst_{\mathcal G_{2}}(n) \leq \St_{\mathcal G_{2}}(n) \leq \mathcal G_{n+2}^{\mathcal L(n)}.
	\]
	Since $\mathcal G$ is a periodic group, so is $\mathcal G_{n+2}$, thus the quotient $\mathcal G_{n+2}/\gamma_3(\mathcal G_{n+2})$ is a finitely generated periodic nilpotent group, hence finite. Consequently $\Rst_{\mathcal G_{2}}(n)$ is of finite index in $\St_{\mathcal G_{2}}(n)$, which is itself of finite index in $\mathcal G_{2}$.
\end{proof}

Note that all branch groups are lawless, cf.\ \cite[Corollary~1.4]{Abe05}, i.e.\ do not satisfy any group identity. This property is clearly bequeathed to supergroups. In particular, the group~$\mathcal{G}$ is neither nilpotent nor an $n$-Engel group for any $n \in \N$.

Now that we have established that the group $\mathcal{G}$ is non-nilpotent and that some of their subgroups are branch groups, it remains to prove that these groups are Engel groups. To achieve this, we want to use \cref{thm:strong contraction > engel}, so we need to establish that $\mathcal{G}$ is an $\mathbf{f}$-Engel group for some $\mathbf{f} \colon \N \to \N_+$, and we need to estimate the function $\mathbf{f}$. To do so, we shall consider certain overgroups of the congruence quotients of $\mathcal{G}$. But first, we need some rather technical lemmata on Engel elements in certain wreath products.

\begin{lemma}\label{lem:abelian wreath}
	Let $r \in \N$ and let $\mathrm{C}_{p^r} = \langle \sigma \rangle$ be a cyclic group of order~$p^r \in \N$ acting faithfully on $X = [1, p^r]$. Let~$G$ be a group and let $W = G \wr_X \mathrm{C}_{p^n}$. Let $g_1, \dots, g_{p^r}, h_1, \dots, h_{p^r} \in G$ be pairwise commuting elements, such that $g_1, \dots, g_{p^r}$ are of order $p$. Then
		\[
			[(g_x)_{x \in X},_{p^r} (h_x)_{x\in X} \sigma] = \id.
		\]
\end{lemma}

\begin{proof}
	For all $n \in \N$ and $x \in X$ the statement
	\[
		([(g, \id, \dots, \id),_n \sigma])\pi_{x} = \prod_{i = 0}^{\lceil n - x/m \rceil} g^{(-1)^{ix} \binom n {ix}},
	\]
	where $\pi_{x}$ denotes the projection to the $x$\textsuperscript{th}~component is easily derived by induction on~$n$, cf.\ \cite[Lemma~3.1]{Lie62}. As a direct application, using well-known divisibility properties of binomial coefficients, we find
	\[
		([(g^{-1}, \id, \dots, \id),_{p^r} \sigma])\pi_{x} = \id
	\]
	for all~$g \in G$ of order $p$. The statement of the lemma can be reduced to the case $(h_x)_{x \in X} = (\id)_{x \in X}$ by observing
	\[
		[(g_x)_{x\in X}, (h_x)_{x \in X}\sigma] = [(g_x)_{x\in X}, \sigma][(g_x)_{x\in X}, (h_x)_{x \in X}]^{\sigma} = [(g_x)_{x\in X}, \sigma]
	\]
	and, using that $[(g_x)_{x \in X}, \sigma]$ is an element of the base group with components commuting with the components of $(h_x)_{x \in X}$,
	\[
		[[(g_x)_{x \in X}, \sigma], (h_x)_{x \in X}] = \id.
	\]
	Also, we see that all elements of the form $(t_x)_{x \in X}$ with $t_x \in \{\id, g_1, \dots, g_{p^r}\}$ for all $x \in X$ commute. Using the easily observed fact that elements of this form are invariant under taking the commutator with $\sigma$, we find
	\begin{align*}
		[(g_x)_{x \in X},_{p^r} \sigma]
		&= [(g_1, \id, \dots, \id),_{p^r} \sigma] [(\id, g_1, \id, \dots, \id,_{p^r} \sigma] \dots [(\id, \dots, \id, g_{p^r}),_{p^r} \sigma] \\&= (\id, \dots, \id),
	\end{align*}
	using that $g_x$ is of order~$p$ for all $x \in X$.
\end{proof}

\begin{lemma}\label{lem:Engel wreath product}
	Let $(A_i)_{i \in \N}$ be a sequence of elementary abelian $p$-groups. Define $W_{0}$ to be the trivial group and $W_{n+1} = W_{n} \wr A_{n}$ for $n \in \N$. The group $W_n$ is a $\frac{p^{n}-1}{p-1}$-Engel group.
\end{lemma}

\begin{proof}
	To (re)gain access to our usual notation we identify~$W_n$ with the group of automorphisms of the finite tree associated to $(A_{n-1}, \dots, A_0, \varnothing, \dots)$ that is generated by the automorphisms $x_{n-1} \dots x_0 \mapsto x_{n-1} \dots x_{i+1} (x_i a_i) x_{i-1} \dots x_0$ for $a_i$ a generator of $A_i$.
	
	We first prove the following claim. Let $g \in \St(k)$ be an element of the $k$\textsuperscript{th} stabiliser and let $h \in W_n$ be any element. Then $[g,_{p^k} h] \in \St(k+1)$. Write $h = (h_x)_{x \in \mathcal{L}(k)} h|^{\mathcal{L}(k)}$. The element $h|^{\mathcal{L}(k)}$ is of order at most~$p^k$ by \cref{lem:max orbit}. Employing~\cref{lem:local checking}, we may restrict to consider $[(g|_x)_{x \in \mathcal{O}}, (h|_x)_{x \in \mathcal{O}}h|^{\mathcal{L}(k)}]$ for all orbits $\mathcal{O}$ of $\mathcal{L}(k)$ under $h|^{\mathcal{L}(k)}$. Calculating modulo $\St(k)'$, this is precisely the situation of \cref{lem:abelian wreath}, hence $[g,_{p^k}h] \in \St(k)'$. But since the action on every layer is abelian, this implies $[g,_{p^k}h] \in \St(k+1)$.
	
	The statement of the lemma now follows easily by induction, since
	\[
		[g,_{1 + p + \dots + p^{n-1}} h] = [\dots[[g,_1 h],_p h] \dots ],_{p^{n-1}} h] \in \St(n+1) = \{\id\}.\qedhere
	\]
\end{proof}

\begin{proposition}\label{prop:engel}
	The group~$\mathcal{G}$ is an Engel group.
\end{proposition}

\begin{proof}
	Define two integer functions $\mathbf{e}, \mathbf{f}\colon \N \to \N$ by $\mathbf{e}(n) = \frac{p^n-1}{p-1}$ and $\mathbf{f}(n) = 2^{\mathbf{e}(n)}$. By \cref{prop:strong contr}, the group~$\mathcal{G}$ is contracting with respect to the function~$\mathbf{g}$ and~$1$. Clearly $\mathbf{f} \prec \mathbf{g}$ (both for even and odd primes) and $\mathbf{e} \prec \mathbf{f}$, thus the group~$\mathcal{G}$ is also contracting with respect to~$\max\{\mathbf{e}, \mathbf{f}\}$ and~$1$. The $n$\textsuperscript{th}~congruence quotient~$\mathcal{G}/\St_{\mathcal{G}}(n)$ is a subgroup of the iterated wreath product $(\dots(X_{n-1} \wr X_{n-2}) \wr X_{n-3}) \dots ) \wr X_1) \wr X_0$, which is a $\mathbf{e}(n)$-Engel group by \cref{lem:Engel wreath product}. Consequently $\mathcal{G}$ is an $\mathbf{e}$-Engel group and we may employ \cref{thm:strong contraction > engel}.
	
	Thus we have to show that there exists some $n \in \N$ (that we shall choose later) such that for all $m \geq n$ and every $h \in \mathcal{G}$ with $\depth(h) \leq m$, every orbit $\mathcal{O}$ of $\mathcal{L}(m)$ under $\tau = h|^{\mathcal{L}(n)}$ and every $g \in \St_{\mathcal{G}}(m)$ with $\depth(g) \leq m$ there exists some $i$ such that
	\[
		[g^{\mathcal{O}},_i h^{\mathcal{O}} \tau] \in \mathcal{G}_m \wr_{\mathcal{O}} \langle \tau \rangle
	\]
	is trivial, where
	\[
		g^{\mathcal{O}} = (g|_x)_{x \in \mathcal{O}} \quad\text{and}\quad h^{\mathcal{O}} = (h|_x)_{x \in \mathcal{O}}.
	\]
	If $[g,_i h] \in \St_{\mathcal{G}}(m+1)$ for some $i \in \N$, then $[g^{\mathcal{O}},_i h^{\mathcal{O}} \tau] \in \St_{\mathcal{G}_m}(1)^{\mathcal{O}}$. Since $\mathcal{G}$ is an $\mathbf{e}$-Engel group,
	\[
		[g^{\mathcal{O}},_{\mathbf{e}(m+1)} h^{\mathcal{O}} \tau] \in \St_{\mathcal{G}_m}(1)^{\mathcal{O}}.
	\]
	Thus we may see this element as a member of $\mathcal{G}_{m+1} \wr_{\{xy \mid x \in \mathcal{O}, y \in X_{m+1}\}} \langle h|^{\mathcal{L}(m+1)} \rangle$, where the action is the one induced by the action of $h$ on $\mathcal{T}$. We again restrict to an orbit: Let $\mathcal{O}^\ast$ be an orbit of $\{xy \mid x \in \mathcal{O}, y \in X_{m+1}\}$ under $\tau^\ast = h|^{\mathcal{L}(m+1)}$ and write
	\[
		g^{\ast} = ([g^{\mathcal{O}},_{\mathbf{e}(m+1)} h^{\mathcal{O}} \tau] \pi_x|_y])_{xy \in \mathcal{O}^\ast} \quad\text{and}\quad h^{\ast} = (h|_{xy})_{xy \in \mathcal{O}^\ast}.
	\]
	By \cref{lem:local checking}, it is enough to prove that the expression
	\[
		[g^\ast,_j h^\ast \tau^\ast]
	\]
	is trivial for some $j \in \N$ to establish that $[g^{\mathcal{O}},_i h^{\mathcal{O}} \tau]$ is trivial for some $i \geq \mathbf{e}(m)$.
	
	First, we notice that, using the same basic technique as in~\cref{lem:upper bound length of iterated word}, the length of $[g^{\mathcal{O}},_{\mathbf{e}(m+1)} h^{\mathcal{O}} \tau] \pi_x$, for any $x \in \mathcal{O}$, is bounded above by
	\[
		2^{\mathbf{e}(m+1)+2} + 2^{\mathbf{e}(m+1)} - 2 \leq 5 \cdot 2^{\mathbf{e}(m+1)} = 5 \cdot \mathbf{f}(m+1).
	\]
	By \cref{lem:max orbit}, the length of $\mathcal{O}^\ast$ is at most $p^{m+1}$. Fix an element $xy \in \mathcal{O}^\ast$ such that $\ell_{E|_{m+1}\cap X_{m+1}}(y)$ is minimal. Since
	\[
		(xy).(h\sigma) = (x.\tau)(y.h|_{x}), \quad (xy).(h\tau)^2 = (x.\tau^2)(y.h_{x}h_{x.\tau}),
	\]
	and so forth, we may describe the elements of $\mathcal{O}^\ast$ by
	\(
		(x.\tau^{k})y_k,
	\)
	where $k \in [0, p^{m+1}-1]$ and 
	\[
		y_k = y.\left(\prod_{j = 0}^{k-1} h_{x.\sigma^j}\right).
	\]
	In particular, $y_0 = y$. We see that
	\begin{equation*}\label{eq:diameter of orbit}
		\dist_{{E|_n} \cap X_n}(y, y_k) = \ell_{E|_{m+1} \cap X_{m+1}}(y^{-1} y_k) = \ell_{E|_n}(y_0^{-1}y_0 \prod_{j = 0}^{k-1} h_{x.\sigma^j}) \leq k \leq p^{m+1},
	\end{equation*}
	using that~$h|_z$ is of word length~$1$ for all $z \in \mathcal{L}(m)$, since $\depth(h) \leq m$.
	
	Now assume that $n \in \N$ is such that for all $m \geq n$ the inequality $\mathbf{d}(m) > 10 \cdot \mathbf{f}(m) + p^m$ holds; this is possible since the function~$\mathbf{d}$ fulfils $\mathbf{f} \prec \mathbf{d}$ and $(n \mapsto p^n) \prec \mathbf{d}$.
	
	We use \cref{lem:separation} with $n = m+1$ and $t = 5 \cdot \mathbf{f}(m+1)$. Since $\dist_{{E|_n} \cap}(y, y_k) \leq p^{m+1}$ and $\mathbf{d}(m+1) - 2t = \mathbf{d}(m+1) - 10 \cdot \mathbf{f}(m+1) > p^{m+1}$, we find that all sections of the form
	\[
		[g^{\mathcal{O}},_{\mathbf{e}(m+1)} h^{\mathcal{O}} \tau] \pi_x|_y \quad\text{or}\quad h|_{xy}
	\]
	for $xy \in \mathcal{O}^\ast$ commute, since the length of $[g^{\mathcal{O}},_{\mathbf{e}(m+1)} h^{\mathcal{O}} \tau] \pi_x$ and of $h|_x$ are bounded by $t$. Thus $g^\ast, h^\ast$ and $\tau^\ast$ meet the conditions of \cref{lem:abelian wreath}, and there exists some $j \in \N$ such that
	\[
		[g^\ast,_j h^\ast\tau^\ast] = \id.
	\]
	This completes our proof.
\end{proof}

We record that \cref{thm:Engel groups} is a direct consequence of \cref{prop:p group}, \cref{prop:branch} and \cref{prop:engel}.


\section{Analogues on $q$-regular trees} 
\label{sec:analogues_on_q_regular_trees}
In this section we prove \cref{thm:Engel p-elements in regular}. The general structure of the proof is the same as for~\cref{thm:Engel groups}; thus we allow ourselves to be more brief in the exposition.

First, we construct the desired groups. Fix an odd prime~$p$. Let~$r \in \N$ be a positive integer, let~$X_r$ be an elementary abelian $p$-group of rank~$r$, generated by $\{ e_1, \dots, e_r \}$, and let $\mathcal T = \mathcal T_{(X_r)_{n \in \N}}$ be the $p^r$-regular rooted tree determined by the constant sequence $(X_r)_{n \in \N}$. Fix the generating set $\widetilde{E} = \bigcup_{i \in [1, r]} \langle e_i \rangle$ for~$X_r$. The diameter of~$X_r$ with respect to~$\widetilde{E}$ is~$r$, and the set~$F$ of elements of furthest $\widetilde{E}$-distance to the identity is $\{ e_1^{k_1} \cdots e_r^{k_r} \mid k_i \in [1, p-1] \text{ for }i \in [1,r] \}$. Thus, viewing~$X_r$ as a vector space over~$\F_p$, there is a basis $V = \{v_1, \dots, v_r\}$ for~$X_r$ contained in $F$. Fix such a basis~$V$.

View the elements of $X_r$ – employing their action on $X_r$ by right-multiplication – as automorphisms of $\mathcal T$ acting on the first letter in a string. Define an automorphism $b \in \St(1)$ by
\[
	b|_x = \begin{cases}
		b &\text{ if }x = \id_{X_r},\\
		e_i &\text{ if }x = v_i \text{ for } i \in [1, r],\\
		\id &\text{ otherwise.}
	\end{cases}
\]
This element is of order $p$. Set $E = \tilde E \cup \langle b \rangle$, and define $\mathcal H = \mathcal H_{p, r} = \langle E \rangle$. The upper companion groups of $\mathcal H$ are easily seen to be equal to $\mathcal H$ itself. For the group $\mathcal{H}$, most (analogues) of the properties we established for $\mathcal{G}$ in the previous section hold true. It is easily seen that $\mathcal H$ is fractal and spherically transitive, using the same techniques as in \cref{lem:basic}. The next two lemmata are proven in the same way as the main claim in the proof of \cref{prop:strong contr} and as \cref{lem:separation}, respectively.

\begin{lemma}
	For all $g \in \mathcal H$ and $u \in \mathcal L(2)$ we have
	\[
		\ell_E(g|_u) \leq \ell_E(g)/r + 1,
	\]
	i.e.\ $\mathcal H$ is contracting with respect to $(n \mapsto r^{n/2})$ and $1$.
\end{lemma}

Just as we have established that $\mathcal G$ is periodic in \cref{prop:p group}, we find that $\mathcal H$ is periodic for $r > p^2$.

\begin{lemma}[Separation]\label{lem:separation 2}
	Let $t \leq r/2$ and let $x \in X_r$. Then
	\[
		\Ball_{\mathcal H}^{E}(t)|_x \subseteq
		\begin{cases}
			\langle b \rangle &\text{ if }\ell_{E}(x) < t,\\
			X_{r} &\text{ if }\ell_{E}(x) > r - t,\\
			\{ \id \} &\text{ otherwise.}
		\end{cases}
	\]
\end{lemma}

As seen in \cref{lem:vanishing commutators}, this implies $[b, e_i, b] = \id$ for $i \in [1, r]$ and $r > 1$.

\begin{proposition}
	The group $\mathcal H_{p, r}$ is a branch group for $r > p^2$ and a weakly branch group for $r > 1$.
\end{proposition}

\begin{proof}
	We proceed similarly as in \cref{prop:branch} and show that $\gamma_3(\upcom_n(\mathcal{H})) = \gamma_3(\mathcal{H})$ is contained in every rigid vertex stabiliser for all $r > 1$. In this case, the rigid layer stabiliser is of finite index if $\mathcal H$ is periodic, i.e.\ if $r > p^2$.
	
	Using the fact that $\mathcal H$ is fractal in the same way as in the proof of \cref{prop:branch}, it is enough to show that for each normal generator $[b, e_i, e_j]$ (recall that $[b, e_i, b]$ is trivial) with $i, j \in [1, r]$ of $\gamma_3(\mathcal{H})$ there exists an element $g \in \St_{\gamma_3(\mathcal{H})}(1)$ such that~$g|_{\id} = [b, e_i, e_j]$. This is a consequence of
	\[
		[b, b^{v_i^{-1}}, b^{v_j^{-1}}] = \begin{cases}
			[b, e_i, e_j] &\text{ if }x = \id,\\
			[b|_{v_i^{-1}}, e_i, b|_{v_i^{-1}v_j}] &\text{ if }x = v_i^{-1},\\
			[b|_{v_j^{-1}}, b|_{v_j^{-1}v_i}, e_j] &\text{ if }x = v_i^{-1},\\
			\id &\text{ otherwise.}
		\end{cases}
	\]
	Since $V$ is linearly independent, all commutator expressions on the right hand side except the first are trivial; one of the entries of each triple commutator is trivial itself.
\end{proof}

Note that one can show, using similar arguments as for the Grigorchuk--Gupta--Sidki groups developed in~\cite{FZ13}, that for sufficiently high $r$ (not depending on $p$) the full commutator subgroup of $\mathcal H_{p,r}$ is contained in every rigid vertex stabiliser; since $\mathcal H_{p,r}$ is finitely generated by finite order elements, the group $\mathcal H_{p,r}$ is branch for all $r$ greater than this uniform bound.

\begin{proposition}\label{prop:e in l}
	Let $r > 3\cdot 2^{p+2}$. Then $\tilde E \subseteq \leftengel(\mathcal H_{p,r})$.
\end{proposition}

\begin{proof}
	Let $g \in \mathcal H_{p,r}$ and let $e \in \tilde E$. We want to prove that there is some integer~$k \in \N$ such that $[g,_k e] = \id$. If $e$ is trivial, we are done. Thus let $e$ have order~$p$. By \cref{lem:Engel wreath product}, $[g,_{p+1} e] \in \St_{\mathcal H_{p,r}}(2)$. Now
	\[
		\ell_{E}([g,_{p+1}, e]|_u) \leq \ell_{E}([g,_{p+1}, e])/r + 1 \leq 2^{p+2}(\ell_{E}(g) + 2)/r + 1 < \ell_{E}(g) + 1
	\]
	for all $u \in \mathcal L(2)$. Applying \cref{lem:local checking}, we may reduce to the case $g = (g_y)_{y \in \mathcal{O}}$ for some orbit $\mathcal{O}$ of $X$ under $e$ with $g_y \in E$ for all $y \in \mathcal{O}$. By \cref{lem:separation 2}, we may suppose that $\{ g_y \mid y \in \mathcal{O} \}$ is contained either in $X_r$ or in $\langle b \rangle$, i.e.\ that there exists $k \in \N$ such that $[g,_k e]$ is contained in $X_r \wr_{\mathcal{O}} \langle e \rangle$ or $\langle b \rangle \wr_{\mathcal{O}} \langle e \rangle$, respectively, which are Engel groups by \cref{lem:abelian wreath}.
\end{proof}

\begin{proposition}\label{prop:b in l}
	Let $r > 3 \cdot 2^{p+2}$. Then $b \in \leftengel(\mathcal H_{p, r})$.
\end{proposition}

\begin{proof}
	We proceed similarly to the last proof. Since $[g, b] \in \St_{\mathcal{H}_{p,r}}(1)$, we have $[g,_n b] = ([[g,b]|_x,_{n-1} b|_x])_{x \in X_r}$ for all $n > 1$. Since all sections of $b$ are either $b$ itself or in $E$, by \cref{prop:e in l} all but the section at $\id$ vanish for sufficiently high $n$. Therefore we only need to consider the section at $\id$. Repeating this argument, we consider $[g,_2 b]|_{\id\id}$. Now
	\[
		\ell_{E}([g, b]|_{\id\id}) \leq \ell_{E}([g, b])/r + 1 \leq 2(\ell_E(g) + 1)/r + 1 < \ell_E(g) + 1,
	\]
	thus we may reduce to $g \in \Ball_{\mathcal H}^{E}(1)$. Since $[e_i, b, b] = \id$ for $i \in [1, r]$, the proof is finished.
\end{proof}

\cref{thm:Engel p-elements in regular} is an immediate consequence of \cref{prop:e in l}, \cref{prop:b in l}, and the fact that the Grigorchuk group is generated by involutions, which are, as demonstrated in~\cite[Theorem~1]{Bar16}, Engel elements. For completion, we quickly give the argument. Given an involution $h$ and any element $g$ of a group $G$, we find
\[
	[g,h]^h = hg^{-1}hg = [h,g],
\]
consequently
\[
	[g,_2 h] = [g,h]^{-1} [g,h]^h = [g,h]^{-2}
\]
and, using induction,
\[
	[g,_{n+1} h] = [g,_n h]^{-1} [g,_nh]^h = [g,h]^{(-2)^{n}}.
\]
If $[g,h]$ is a $2$-element, the involution $h$ is left Engel on $g$. Thus, in a $2$-group -- such as the Grigorchuk group -- every involution is a left Engel element.


\section{Further remarks and questions} 
\label{sec:further_remarks_and_questions}

\subsection*{Comparison to the examples of Golod} 
\label{sub:comparison_to_the_examples_of_golod}

In his 1969 paper \cite{Gol69}, Golod constructs a family $(G_d)_{d \in \N_{>2}}$ of finitely generated non-nilpotent Engel groups. Since this family and the groups constructed in \cref{sec:construction_of_finitely_generated_non_nilpotent_engel_groups} are the only known examples with these properties, a comparison seems appropriate.

The groups constructed by Golod have the following property. The non-nilpotent group $G_d$ is $d$-generated, and every $(d-1)$-generated subgroup is nilpotent; in this way, his counterexample (for $d > 3$) is more than `just' an Engel group. The groups $\mathcal{G}_p$ do not have this property: Indeed, they contain non-nilpotent subgroups of arbitrarily large rank, namely their upper companion groups. In this sense, our groups are (heuristically) `less nilpotent'.

Golod claims that the methods developed in~\cite{Gol69} suffice to construct a $2$-generated non-nilpotent Engel group, although he gives no proof. A benefit of our construction is that the group $\mathcal{G}_p$, at least for odd primes, is given by a explicitly defined generating set consisting of two elements of order~$p$.

It would be interesting to quantify how `strongly' both families of groups (or any other non-nilpotent Engel groups) satisfy the Engel property. We suggest the following definition. Let~$G$ be an Engel group generated by a finite set~$S$. For every $n \in \N$, let $\mathbf{e}_G(n)$ be the least integer such that
\[
	[g,_{\mathbf{e}(n)} h] = \id
\]
for all $g, h \in \Ball_G^S(n)$. The \emph{Engel growth function} $\mathbf{e}_G\colon \N \to \N$ is bounded by $c \in \N$ if $G$ is nilpotent of class~$c$. In general, the growth of the function $\mathbf{e}_G$ gives a quantification of the Engel property, as an analogue to the periodicity growth function introduced by Grigorchuk in~\cite{Gri83}, that gives a measure on how periodic a group~$G$ is.

\begin{question}
	How do the growth types of the functions $\mathbf{e}_{\mathcal{G}_p}$ and $\mathbf{e}_{G_d}$ compare?
\end{question}


\subsection*{Engel groups acting on $q$-regular trees} 
\label{sub:engel_groups_q_regular_trees}

Our methods cannot provide an Engel group on a $q$-regular tree, i.e.\ \cref{thm:Engel p-elements in regular} cannot be easily improved. To see this, let~$G$ an Engel finitely generated branch group contracting with respect to a super-exponentially growing function~$\mathbf{g}\colon \N \to \N$. Then for every $\mu > 0$ there exists a sufficiently high integer~$n \in \N$ such that
\[
	\sum_{u \in \mathcal{L}(n)} \ell_{E}(g|_u) \leq |\mathcal{L}(n)| \cdot (\ell_{E}(g)/(\mathbf{g}(n)) + 1 )< \mu \ell_E(g) + |\mathcal{L}(n)|
\]
for all $g \in G$. Using Grigorchuk's `strong contraction' argument to establish intermediate growth for subgroups of $\Aut(\mathcal T)$, see \cite{Gri83} or \cite{BGS03}, we find that such a group has a word growth function that is bound from above by $e^{n^\mu}$ for arbitrarily small $\mu$.

Since every congruence quotient of~$G$ is a finite Engel group, hence nilpotent, the group~$G$ is residually nilpotent. But by the solution of the gap conjecture for finitely generated residually nilpotent groups, cf.\ \cite[Theorem~10.2]{Gri14} or \cite{Wil11} for a more general result, every such group with super-polynomial word growth has growth rate at least $e^{\sqrt n}$. Since branch groups contain direct products of a subgroup of finite index as a subgroup of finite index, they cannot be of polynomial growth. Thus, such a group~$G$ cannot exist.

However, this does not exclude the possibility that the groups $\mathcal H_{p,r}$ are Engel. Therefore, the following question is pressing:
\begin{question}
	Does there exist a pair of a prime~$p$ and an integer~$r$ such that $\mathcal H_{p,r}$ is an Engel group?
\end{question}
Or, more generally:
\begin{question}
	Does there exist a finitely generated Engel branch group acting on a regular rooted tree?
\end{question}
Note that if $\mathcal H_{p,r}$ is not an Engel group for some sufficiently high $r$, it is an non-Engel group generated by Engel elements. All known examples of such groups are $2$-groups; it is an open question due to Bludov if there exist non-$2$-groups with this property.

\subsection*{Ubiquity of Engel branch groups} 
\label{sub:ubiquity_of_engel_branch_groups}

The construction exhibited in \cref{sec:construction_of_finitely_generated_non_nilpotent_engel_groups} can clearly be varied to produce more Engel branch groups, for example by replacing the sets~$F(n)$ by sets with similar properties, or (less trivially) by replacing the tree build out of elementary abelian $p$-groups with a tree build out of other finite Engel groups. In view of the classic result of Gupta and Sidki~\cite{GS83a} that every finite $p$-group is contained in the corresponding Gupta--Sidki $p$-group (a phenomenon that is not uncommon among branch groups), we ask the following question.
\begin{question}
	Does there exist a finitely generated Engel branch group containing all finite Engel groups as a subgroup?
\end{question}
Due to a recent result of Kionke and Schesler \cite{KS22}, every finitely generated residually finite periodic group embeds into a finitely generated perfect periodic branch groups. Thus we ask:
\begin{question}
	Does every finitely generated residually finite Engel group embed into a finitely generated Engel branch group?
\end{question}


\subsection*{Other iterated identities and their hierarchy} 
\label{sub:other_iterated_identities_and_their_hierarchy}

In view of \cref{thm:strong contraction > engel}, it is natural to ask the following question:

\begin{question}
	Let $w \in F_{k+1}$ be a word. Is there a finitely generated Engel-$w$ branch group?
\end{question}

Naturally, the problem of constructing such groups is related to finding bounds on $n$ for which iterated wreath products are $n$-Engel-$w$ groups. 

Note that Fernández-Alcober, Noce and Tracey \cite[Theorem~B]{FNT20} proved that every Engel branch group is periodic. While our methods do not yield this result, they help to explain it heuristically: the contraction needed to be periodic is weaker than the contraction necessary for being Engel. Motivated by their result, we ask:

\begin{question}
	Let $G$ be an Engel-$w$ branch group for some word $w \in F_{k+1}$. Is $G$ always periodic? Are there words~$w$ other than the commutator which enforce that $G$ is an Engel group?
\end{question}



\textbf{Declarations of interest}: none.



\begin{thebibliography}{10}

\bibitem{Abe05}
Mikl{\'o}s Ab{\'e}rt, \emph{Group laws and free subgroups in topological
  groups}, Bulletin of the London Mathematical Society \textbf{37} (2005),
  no.~04, 525--534.

\bibitem{Bar16}
Laurent Bartholdi, \emph{Algorithmic {{Decidability}} of {{Engel}}'s
  {{Property}} for {{Automaton Groups}}}, Computer {{Science}} \textendash{}
  {{Theory}} and {{Applications}} (Alexander~S. Kulikov and Gerhard~J.
  Woeginger, eds.), vol. 9691, {Springer International Publishing}, {Cham},
  2016, pp.~29--40.

\bibitem{BGS03}
Laurent Bartholdi, Rostislav~I. Grigorchuk, and Zoran {\v S}uni{\'k},
  \emph{Branch groups}, Handbook of {{Algebra}}, vol.~3, {Elsevier}, 2003,
  pp.~989--1112.

\bibitem{Blu05}
V.~V. Bludov, \emph{An example of not {{Engel}} group generated by {{Engel}}
  elements}, A {{Conference}} in {{Honor}} of {{Adalbert Bovdi}}'s 70th
  {{Birthday}} ({Debrecen, Hungary}), November 2005, pp.~7--8.

\bibitem{Ers15}
Anna Erschler, \emph{Iterated identities and iterational depth of groups},
  Journal of Modern Dynamics \textbf{9} (2015), no.~01, 257--284.

\bibitem{FZ13}
Gustavo {Fern{\'a}ndez-Alcober} and Amaia {Zugadi-Reizabal},
  \emph{{{GGS-groups}}: {{Order}} of congruence quotients and {{Hausdorff}}
  dimension}, Transactions of the American Mathematical Society \textbf{366}
  (2013), no.~4, 1993--2017.

\bibitem{FGN19}
Gustavo~A. {Fern{\'a}ndez-Alcober}, Albert Garreta, and Marialaura Noce,
  \emph{Engel elements in some fractal groups}, Monatshefte f\"ur Mathematik
  \textbf{189} (2019), no.~4, 651--660.

\bibitem{FNT20}
Gustavo~A. {Fern{\'a}ndez-Alcober}, Marialaura Noce, and Gareth~M. Tracey,
  \emph{Engel elements in weakly branch groups}, Journal of Algebra
  \textbf{554} (2020), 54--77.

\bibitem{Gol65}
E.~S. Golod, \emph{On nil-algebras and finitely approximable
  {${p}$}-groups}, vol.~48, pp.~103--106, {American Mathematical
  Society}, {Providence, Rhode Island}, 1965.

\bibitem{Gol69}
\bysame, \emph{Some problems of {{Burnside}} type}, American {{Mathematical
  Society Translations}}: {{Series}} 2, vol.~84, {American Mathematical
  Society}, {Providence, Rhode Island}, 1969, pp.~83--88.

\bibitem{Gri83}
Rostislav Grigorchuk, \emph{On the {{Milnor}} problem of group growth}, Dokl.
  Akad. Nauk SSSR \textbf{271} (1983), no.~1, 30--33.

\bibitem{Gri14}
\bysame, \emph{Milnor's problem on the growth of groups and its consequences},
  pp.~705--774, {Princeton University Press}, December 2014.

\bibitem{GS83a}
Narain Gupta and Said Sidki, \emph{{Some infinite p-groups}}, Algebra i Logika
  \textbf{22} (1983), no.~5, 584--589.

\bibitem{KS22}
Steffen Kionke and Eduard Schesler, \emph{Realizing residually finite groups as
  subgroups of branch groups}, December 2022.

\bibitem{Lie62}
Hans Liebeck, \emph{Concerning nilpotent wreath products}, Mathematical
  Proceedings of the Cambridge Philosophical Society \textbf{58} (1962), no.~3,
  443--451.

\bibitem{Nek05}
Volodymyr Nekrashevych, \emph{Self-{{Similar Groups}}}, Mathematical
  {{Surveys}} and {{Monographs}}, vol. 117, {American Mathematical Society},
  {Providence, Rhode Island}, August 2005.

\bibitem{NT18}
Marialaura Noce and Antonio Tortora, \emph{A note on {{Engel}} elements in the
  first {{Grigorchuk}} group}, International Journal of Group Theory (2018),
  no.~Online First.

\bibitem{Pet23}
J.~Moritz Petschick, \emph{Two periodicity conditions for spinal groups},
  Journal of Algebra \textbf{633} (2023), 242--269.

\bibitem{Sid87}
Said Sidki, \emph{On a 2-generated infinite 3-group: {{Subgroups}} and
  automorphisms}, Journal of Algebra \textbf{110} (1987), no.~1, 24--55.

\bibitem{TT20}
Antonio Tortora and Maria Tota, \emph{Engel groups in bath - ten years later},
  International Journal of Group Theory \textbf{9} (2020), no.~4.

\bibitem{Tra11}
Gunnar Traustason, \emph{Engel groups}, Groups {{St}}. {{Andrews}} 2009.
  {{Vol}}. {{II}}. {{Selected}} Papers of the Conference, {{University}} of
  {{Bath}}, {{Bath}}, {{UK}}, {{August}} 1--15, 2009., {Cambridge: Cambridge
  University Press}, 2011, pp.~520--550.

\bibitem{Wil11}
John~S. Wilson, \emph{The gap in the growth of residually soluble groups},
  Bulletin of the London Mathematical Society \textbf{43} (2011), no.~3,
  576--582.

\end{thebibliography}
\end{document}